\nonstopmode
\documentclass[a4paper, 10pt]{amsart}
\usepackage{latexsym}
\usepackage{amsmath, amssymb}
\usepackage[ansinew]{inputenc}
\usepackage[all]{xy}
\theoremstyle{plain}
\newtheorem*{definition*}{Definition}

\newtheorem*{lemma*}{Lemma}
\newtheorem{lemma}[subsection]{Lemma}
\newtheorem*{theorem*}{Theorem}
\newtheorem{theorem}[subsection]{Theorem}
\newtheorem*{proposition*}{Proposition}
\newtheorem{proposition}[subsection]{Proposition}
\newtheorem*{corollary*}{Corollary}
\newtheorem{corollary}[subsection]{Corollary}
\theoremstyle{remark}
\newtheorem*{remark*}{Remark}
\newtheorem{remark}[subsection]{Remark}
\theoremstyle{remark}
\newtheorem*{example*}{Example}
\newtheorem{example}[subsection]{Example}
\def\o{\circ}
\def\al{\alpha}

\def\ga{\gamma}
\def\de{\delta}

\def\vartheta{\theta}

\def\ph{\varphi}

\def\ps{\psi}
\def\om{\omega}

\def\Om{\Omega}

\def\x{\times}
\def\g{\mathfrak g}

\def\on{\operatorname}

\title[Cohomology for a group
of diffeomorphisms 
]
{Cohomology for a group
of diffeomorphisms of a manifold preserving an exact form.}
\author[Losik, Michor]{Mark Losik, Peter W. Michor}

\thanks{Supported
by FWF Projekt P~17108-N04.}
\address{M. Losik: Saratov State University, Astrakhanskaya, 83,
410026 Saratov, Russia.}
\email{LosikMV@info.sgu.ru}
\address{P.\  W.\  Michor: Fakult\"at f\"ur Mathematik, Universit\"at Wien,
Nordbergstrasse 15, A-1090 Wien, Austria; {\it and}:
Erwin Schr\"odinger Institute of Mathematical Physics, Boltzmanngasse
9, A-1090 Wien, Austria.}
\email{Peter.Michor@esi.ac.at}
\date{\today}

\begin{document}
\begin{abstract}
Let $M$ be a $G$-manifold and $\om$ a $G$-invariant exact $m$-form on $M$. 
We indicate when these data allow us to constract a cocycle on a group $G$ with 
values in the trivial $G$-module $\mathbb R$ and when this cocycle is nontrivial. 
\end{abstract}
\maketitle

\section{Introduction} Let $M$ be a manifold, let $G$ be a group of diffeomorphisms 
of $M$, and let $\om$ be a $G$-invariant exact $m$-form on $M$. In this paper we apply 
the construction of \cite{Lo} to get from these data a cocycle on 
the group $G$ with values in the trivial $G$-module $\mathbb R$. 
We prove that this cocycle may be chosen
differentiable (continuous) whenever $G$ is a subgroup of a Lie group 
(a topological group). Moreover, we prove that for    
a manifold $\mathbb R^n\x M$ with an exact form $\om$ which is either of
type $\om_0+\om_M$ or of type $\om_0\wedge\om_M$, where $\om_0$ is a
nonzero form on $\mathbb R^n$ with constant coefficients and $\om_M$ is a
form on $M$, and the group $\on{Diff}(\mathbb R^n\x M,\om)$ of diffeomorphisms
of $\mathbb R^n\x M$ preserving the form $\om$ the corresponding cocycle is nontrivial.

\section{A construction of cohomology classes for a group of diffeomorphisms} 

Let $G$ be a group and let $A$ be a right $G$-module. Let
$C^p(G,A)$ be the set of maps from $G^p$ to $A$ for $p>0$ and let
$C^0(G,A)=A$. Define the differential $D:C^p(G,A)\to C^{p+1}(G,A)$ as
follows, for $f\in C^p(G,A)$ and $g_1,\dots,g_{p+1}\in G$:
\begin{multline}\label{de'}
(Df)(g_1,\dots,g_{p+1})=f(g_2,\dots,g_{p+1})\\+
\sum_{i=1}^p(-1)^if(g_1,\dots,g_{i}g_{i+1},\dots,g_{p+1})+
(-1)^{p+1}f(g_1,\dots,g_p)g_{p+1}.
\end{multline}
Then $C^*(G,A)=(C^p(G,A),D)_{p\ge 0}$ is the standard complex of
nonhomogeneous cochains of the group $G$ with values in the right
$G$-module $A$ and its
cohomology $H^*(G,A)=(H^p(G,A)_{p\ge 0})$ is the cohomology of the group $G$
with values in $A$.

Let $M$ be a smooth $n$-dimensional $G$-manifold, where $G$ is a group of
diffeomorphisms of $M$. Denote by $\Om(M)=(\Om^p(M))_{p=1,\dots,n}$ the 
de~Rham complex of differential forms on $M$ and consider the natural (right)
action of the group $G$ on $\Om(M)$ by pull backs. Denote by $\Om(M)^G$ the
subcomplex of $\Om(M)$ consisting of $G$-invariant forms.
Next we denote by $H_q(M)$ the $q$-dimensional real homology of $M$ and by
$H^q(M)$ the $q$-dimensional real cohomology of $M$.

Let $C^*(G,\Om(M))=\{C^p(G,\Om^q(M)),\delta')\}_{p,q\ge 0}$ be 
the standard complex of nonhomogeneous cochains of $G$ with values in
the $G$-module $\Om(M)$. 
We define the second differential
$\de'':C^p(G,\Om^q(M))\to C^p(G,\Om^{q+1}(M))$ by 
$$
(\de''c)(g_1,\dots,g_p)=(-1)^pdc(g_1,\dots,g_p),
$$
where $f\in C^p(G,\Om^q(M))$, $g_1,\dots,g_p\in G$, and where $d$ is exterior
derivative. Since $\de'\de''+\de''\de'=0$,
we have on $C^*(G,\Om(M))$ the structure of double complex. 
Denote by $C^{**}(G,\Om(M))$ the cochain complex $C^*(G,\Om(M))$ 
with respect to the total differential $\de=\de'+\de''$.
We denote by $H(G,M,\Om(M))$ the cohomology of the complex
$C^{**}(G,\Om(M))$.

It is easily checked that the inclusion 
$\Om(M)^G\subset C^0(G,\Om(M))$ induces an injective homomorphism of complexes
$\Om(M)^G\to C^{**}(G,\Om(M))$ and thus also a homomorphism 
$H(\Om(M)^G)\to H(G,M,\Om(M))$ of cohomologies. 
We identify $\om\in\Om(M)^G$ with its image by the inclusion 
$\Om(M)^G\subset C^{**}(G,\Om(M))$ and denote
by $h(\om)$ the cohomology class of $\om$ in the complex $C^{**}(G,\Om(M))$ 
whenever the form $\om$ is closed.
 
We shall use some standard facts on spectral sequences 
(see, for example, \cite{C-Ei}). 
Consider the first filtration 
$$
F_{1,p}(G,M,\Om(M)):=\oplus_{q\ge p}C^q(G,\Om(M))
$$ 
of the double complex $C^*(G,\Om(M))$. 
By definition,  
$F_{1,p}(G,M,\Om(M))$ is a subcomplex of the complex $C^{**}(G,\Om(M))$ and
$F_{1,0}(G,M,\Om(M))=C^{**}(G,\Om(M))$. Denote by $E_{1,r}=(E_{1,r}^{p,q},
d_{1,r}^{pq})_{p,q\ge 0}$ for $r=0,1,\dots,\infty$ the corresponding
spectral sequence.
Denote by $h_p$ the homomorphism of cohomologies $H(F_{1,p}(G,M,\Om(M))\to
H(G,M,\Om(M))$ induced by the inclusion $F_{1,p}(G,M,\Om(M))\subset
C^{**}(G,\Om(M))$. Then $(\on{Im}h_p)_{p\ge 0}$ is a filtration of the
cohomology $H(G,M,\Om(M))$ and
$$
E^{p,q}_{1,\infty}=
h_p(H^{p+q}(F_{1,p}(G,M,\Om(M))))/h_{p+1}(H^{p+q}(F_{1,p+1}(G,M,\Om(M)))).
$$
For this spectral sequence we have $E^{p,q}_{1,2}=H^p(G,H^q(M))$, where
an action of the group $G$ on $H^q(M)$ is induced by its action on
$\Om(M)$. 
Moreover, there is a natural homomorphism 
$H^p(G,H^0(M))=E^{p,0}_{1,2}\to E^{p,0}_\infty$. 

\begin{proposition}\label{m=p+1} 
Let $\om\in\Om(M)^G$ be an exact $m$-form and let $p$ be the maximal integer
such that $h(\om)\in h_{p+1}(H^m(F_{1,p+1}(G,M,\Om(M))))$.
Then the image of $h(\om)$ under the natural homomorphism 
$\on{Im}h_{p+1}\to\on{Im}h_{p+1}/\on{Im}h_{p+2}$ 
belongs to $E^{p+1,m-p-1}_{1,m-p+1}$. In particular, if  
$m=p+1$ and the manifold $M$ is connected, the above image of $h(\om)$ 
is an $m$-dimensional cohomology class of the group $G$ with values in the
trivial $G$-module $\mathbb R$.
\end{proposition}

\begin{proof} By assumption the image of $h(\om)$ under the homomorphism 
$$
\on{Im}h_{p+1}\to\on{Im}h_{p+1}/\on{Im}h_{p+2}
$$ 
belongs to $E^{p+1,m-p-1}_{1,\infty}$.
Since $d_{1,r}^{p,q}:E_{1,r}^{p,q}\to E_{1,r}^{p+r,q-r+1}$ vanishes
whenever $r>q+1$, we
have $E^{p+1,m-p-1}_{1,\infty}=E^{p+1,m-p-1}_{1,m-p+1}$. 
If $m=p+1$  and the manifold $M$ is connected we have 
$E^{p+1,m-p-1}_{1,\infty}=E^{m,0}_2=H^m(G,H^0(M))=H^m(G,\mathbb R)$. 
\end{proof}

\begin{theorem}\label{group cocycle} 
Let $\om$ be a $G$-invariant exact $m$-form on $M$ and let
$H^{m-p}(M)=\dots =H^{m-1}(M)=0$ and
$H^{m-p-1}(M)\ne 0$ for some $1\le p\le m-1$. Then
$h(\om)\in\on{Im}h_{p+1}$
and $h(\om)$ defines a unique $(p+1)$-dimensional cohomology class $c(\om)$
of the group $G$ with values in the natural $G$-module $H^{m-p-1}(M)$.
\end{theorem}
\begin{proof}
By assumption there is a form 
$\ph_{0,m-1}\in C^0(G,\Om^{m-1}(M))$ such that
$\om=-d\ph_{0,m-1}=-\de''\ph_{0,m-1}$. Then we have
$\om+\de\ph_{0,m-1}=\de'\ph_{0,m-1}$.

Since $H^{m-1}(M)=0$ and 
$\de''\de'\ph_{0,m-1}=-\de'\de''\ph_{0,m-1}=\de'\om=0$,
there is a cochain $\ph_{1,m-2}\in C^1(G,\Om^{m-2}(M))$ such that 
$\de'\ph_{0,m-1}=-\de''\ph_{1,m-2}$. Thus we have
$\om+\de(\ph_{0,m-1}+\ph_{1,m-2})=\de'\ph_{1,m-2}$.

Using the conditions 
$H^{m-2}(M)=\dots =H^{m-p}(M)=0$ and proceeding in the same way we get 
for $i=1,\dots,p$  the cochains $\ph_{i,m-i-1}\in C^i(G,\Om^{m-i-1}(M))$ 
such that
\begin{equation}\label{sequence}
\de'\ph_{i-1,m-i}+\de''\ph_{i,m-i-1}=0
\end{equation}
and so  
$$
\om+\de(\ph_{0,m-1}+\dots+\ph_{p,m-p-1})=\de'\ph_{p,m-p-1}\in C^{p+1}(G,\Om^{m-p-1}(M)).
$$ 
Moreover, we have 
$$
d\de'\ph_{p,m-p-1}=-\de''\de'\ph_{p,m-p-1}=\de''\de''\ph_{p-1,m-p}=0.
$$ 
Consider $H^{m-p-1}(M)$ as a $G$-module with respect to the natural action
of $G$ on $H^{m-p-1}(M)$.
Then the cochain $\de'\ph_{p,m-p-1}$ defines a cocycle on $G$ of degree
$p+1$ with values
in the $G$-module $H^{m-p-1}(M)$.  
We claim that the cohomology class of this cocycle depends only on the
cohomology class of $\om$ in the complex $\Om(M)^G$.

If we replace the form $\om$ by a form $\om+d\om_1$, where 
$\om_1\in\Om^{m-1}(M)\cap\Om(M)^G$, one can replace the sequence
$\ph_{0,m-1},\dots,\ph_{p,m-p-1}$ by the sequence
$\ph_{0,m-1}-\om_1,\ph_{1,m-2},\dots,\ph_{p,m-p-1}$ and obtain the same cochain
$\ph_{p,m-p-1}$ at the end. 

Consider another sequence $\bar\ph_{0,m-1},\dots,\bar\ph_{p,m-p-1}$
$(i=0,\dots,p)$ such that $\om=-d\bar\ph_{0,m-1}$ and
$\de'\bar\ph_{i-1,m-i}+\de''\bar\ph_{i,m-i-1}=0$ for $i=1,\dots,p$.
Since $H^{m-1}(M)=0$ we have 
$$
\bar\ph_{0,m-1}=\ph_{0,m-1}+\de''\ps_{0,m-2},
$$ 
where $\ps_{0,m-2}\in C^0(G,\Om^{m-2}(M))$. If $p=1$ we have 
$\de'\bar \ph_{0,m-1}=\de'\ph_{0,m-1}-\de''\de'\ps_{0,m-2}$ and we are done. 
If $p>1$ we have  
$$
\de'\bar\ph_{0,m-1}=\de'\ph_{0,m-1}-\de''\de'\ps_{0,m-2}=-\de''(\ph_{1,m-2}
+\de'\ps_{0,m-2})=-\de''\bar\ph_{1,m-2}.
$$
Since $H^{m-2}(M)=0$ there is a cochain $\ps_{1,m-3}\in
C^1(G,\Om^{m-3}(M))$ such that
$\bar\ph_{1,m-2}=\ph_{1,p-2}+\de'\ps_{0,m-2}+\de''\ps_{1,m-3}$.
For $i=2,\dots,p-1$ proceeding in the same way we get the cochains
$\ps_{i,m-i-2}\in C^i(G,\Om^{m-i-2}(M))$ such that
$$
\bar\ph_{i,m-i-1}=\ph_{i,m-i-1}+\de'\ps_{i-1,m-i-1}+\de''\ps_{i,m-i-2}.
$$ 
In particular, we have 
$$
\bar\ph_{p,m-p-1}=\ph_{p,m-p-1}+\de'\ps_{p-1,m-p-1}+\de''\ps_{p,m-p-2}
$$    
and $\de'\bar\ph_{p,m-p-1}=\de'\ph_{p,m-p-1}-\de''\ps_{p,m-p-2}$. Thus the
cochains $\de'\bar\ph_{p,m-p-1}$ and $\de'\ph_{p,m-p-1}$ define the same
cohomology class of $H^{p+1}(G,H^{m-p-1}(M))$.
\end{proof}

Suppose that the conditions of theorem \ref{group cocycle} are satisfied
for an exact $m$-form $\om\in\Om(M)^G$.  
Let $\al$ be a singular smooth cycle of $M$ of dimension $m-p-1$ whose
homology class $a$ is invariant under the natural action of the group $G$
on $H_{m-p-1}(M)$.
Put 
\begin{equation}\label{c_a}
c_a(\om)(g_1,\dots,g_{p+1})=\int_{\al}(\de'\ph_{p,m-p-1})(g_1,\dots,g_{p+1}).
\end{equation}
By definition, $c_a(\om)$ is a $(p+1)$-cocycle on the group $G$ with values
in the trivial $G$-module $\mathbb R$ which is independent of a choice of
the cycle $\al$ in the homology
class $a$. 

Let $p=0$. Evidently, the cocycle $c_a(\om)$ is nontrivial if and only if it
does not vanish. From now on we assume $p>0$.

\begin{remark} 
Let the assumptions of theorem \ref{group cocycle} be satisfied for an 
exact $m$-form $\om\in\Om(M)^G$. If either the manifold is connected and
$m=p+1$, or $G$ is a connected topological group, the action of $G$ on the
$H_{m-p-1}(M)$ is trivial.
If the homology class $a$ of the cycle $\al$ is not invariant under  the
action of $G$, consider the vector subspace $H_a$ of $H_{m-p-1}(M)$
generated by the orbit of $a$. Then (\ref{c_a}) defines a $(p+1)$-cocycle
on the group $G$ with values in the $G$-module $H_a$.
\end{remark}

Consider the partial case when $M=G$ is a connected Lie group and the group
$G$ acts on
$M$ by left translations. It is clear that the complex $\Om(G)^G$ is
isomorphic to the complex $C^*(\g,\mathbb R)$ of standard cochains of the
Lie algebra $\g$ of the group $G$ with values in the trivial $\g$-module
$\mathbb R$. Consider the second filtration 
$$
F_{2,p}C^{**}(G,G,\Om(G)):=\oplus_{q\ge p} C^*(G,\Om^q(G))
$$ 
of the double complex $C^{**}(G,\Om(G))$   
and the corresponding spectral sequence
$E_{2,r}=(E^{p,q}_{2,r},d^{p,q}_r)_{p,q\ge 0}$ for $r=0,1,\dots,\infty$. It
is easily seen that $E_{2,1}^{p,q}=H^p(G,\Om^q(G))$.

\begin{lemma}\label{H(G,Om(G))} 
The inclusion $\Om(M)^G\subset C^{**}(G,G,\Om(G))$ induces an isomorphism of 
cohomologies.
\end{lemma} 
\begin{proof} We prove that for each $q\ge 0$ we have 
$$
H^p(G,\Om^q(G))=0\quad\mbox{for $p>0$ and $H^0(G,\Om^q(G))=\Om^q(G)^G$}.
$$

First we consider the case when $q=0$. 
We use the standard operator $B:C^p(G,\Om^0(G))\to C^{p-1}(G,\Om^0(G))$
defined as follows. For $p>0$, put
$$
(Bc)(g_1,\dots,g_{p-1})(g)=(-1)^pc(g_1,\dots,g_{p-1},g)(e),
$$ 
where $c\in C^p(G,\Om^0(G))$, $g,g_1,\dots,g_{p-1}\in G$, and $e$ is the
identity element of $G$. For $c\in C^0(G,\Om^0(G))$, put $Bc=0$. It is easy
to check that $B$ is a homotopy operator between the identity isomorphism
of $C^*(G,\Om^0(G))$ and the map
of $C^*(G,\Om^0(G))$ into itself which vanishes on $C^p(G,\Om^0(G))$ for
$p>0$ and takes $c\in C^0(G,\Om^0))$ to $c(e)$. 
This proves our statement for $p=0$.

To prove our statement for $p>0$ we note that
$\Om^q(G)=\Om^q(G)^G\otimes\Om^0(G)$, where
$\Om^p(G)^G$ is the space of left invariant $q$-forms on $G$. Since $G$
acts trivially
on $\Om^p(G)^G$, its action on $\Om^q(G)$ is induced by its action on
$\Om^0(G)$. Then we have $H^p(G,\Om^q(G))=0$ for $p>0$ and
$H^0(G,\Om^q(G))=\Om^q(G)^G$.

The above statement implies that $E^{p,q}_{2,1}=0$ when $p>0$ and
$E^{0,q}_{2,1}=\Om^q(G)^G$. Then $E_{2,1}=\Om(G)^G$ and evidently the
differential $d^{0,q}_{2,1}$ equals the exterior derivative $d$ on
$\Om(G)^G$ up to sign. Therefore we have $E_{2,2}^{p,q}=0$ when $p>0$ and
$E^{0,q}_{2,2}=H^q(\Om(G)^G)$. Thus implies that
$E_{2,\infty}^{p,q}=E_{2,2}^{p,q}$ and therefore
the inclusion $\Om(G)^G\subset C^{**}(G,\Om(G))$ induces an isomorphism of
cohomologies.
\end{proof}

\begin{proposition}\label{nontriviality} 
Let $\om\in\Om(G)^G$ be an exact $m$-form whose cohomology class in the
complex $\Om(G)^G$ is nontrivial and let
$$
\ph_{i,m-i-1}\in C^i(G,\Om^{m-i-1}(G))\quad(i=0,\dots,m-1)
$$ 
be a sequence of cochains such that 
$\de(\om+\ph_{0,m-1}+\dots+\ph_{m-1,0})=\de'\ph_{m-1,0}$. 
Then, for a point $x\in G$, a cocycle $\de'\ph_{m-1,0}(g)(x)$ of the complex 
$C^*(G,\mathbb R)$ is nontrivial.
\end{proposition}

\begin{proof} By lemma \ref{H(G,Om(G))} the cohomology class of $\om$ in
the complex
$C^{**}(G,\Om(G))$ is nontrivial and then by assumption $\om$ defines a
nontrivial element of $E_{2,\infty}^{m,0}$. Since
$$
H^m(G,\mathbb R)=E_{2,2}^{m,0}=E_{2,\infty}^{m,0},
$$
the cocycle $\de'\ph_{m-1,0}(g)(x)$ of the complex 
$C^*(G,\mathbb R)$ is nontrivial.  
\end{proof}  

\section{The map $f_\ga$ and its properties}
 
Let $G$ be a finite dimensional Lie group. 
For $X\in T_e(M)$ denote by $X^r$ the right invariant vector field on $G$
such that
$X^r(e)=X$ and by $\tilde X$ the fundamental vector field on $M$
corresponding to $X$
for the action of $G$ on $M$. 
We denote the action by $\ph:G\x M\to M$ and write
$gx=\ph(g,x)=\ph^x(g)=\ph_g(x)$.
By definition, $T(\ph^x)X^r(g)=\tilde X(gx)$ and for each $g\in
G$ we have
$\ph^x\o \on{L}_g=\ph_g\o\ph^x$, where $\on{L}_g$ is left translation on $G$. 

Let $\ga$ be a singular smooth cycle of dimension $q$ on $M$. Define a map
$f_\ga:\Om(M)\to\Om(G)$
as follows. Let $\om\in\Om(M)$. If $\deg\om<q$ put $f_\ga(\om)=0$.  
If $\deg\om=p+q$ with $p\ge 0$ put 
\begin{equation}\label{f}
f_\ga(\om)(X_1^r,\dots,X_p^r)(g)=\int_{\ga}\ph_g^*\left(i(\tilde X_p)\dots 
i(\tilde X_1)\om\right) =
\int_{g\ga}i(\tilde X_p)\dots i(\tilde X_1)\om,
\end{equation}
where $X_1,\dots, X_p\in T_e(G)$ and $g\in G$. Clearly $\om\to f_\ga(\om)$
is a linear map from $\Om(M)$ to $\Om(G)$ decreasing degrees to $q$.

Consider the action of the group $G$ on itself by left translations. 

\begin{lemma}\label{equivariant} 
The map $f_\ga$ is $G$-equivariant.
\end{lemma}

\begin{proof} 
It suffices to consider the case when $\deg\om\ge q$. Let
$\om\in\Om^{p+q}(M)$, $X\in T_e(G)$, and $g,\tilde g\in G$.
It is easy to check that 
\begin{align}
T(\on{L}_{\tilde g})\o X^r
&=(\on{ad}\tilde g(X))^r\o L_{\tilde g}:G\to TG,
\label{1}\\
T(\ph_{\tilde g})\o \tilde X 
&= \widetilde{\on{ad}\tilde g(X)}\o \ph_{\tilde g}:M\to TM
\label{2}\end{align}
From this we get
\begin{multline*}
L_{\tilde g}^*f_\ga(\om))(X_1^r,\dots,X_p^r)(g)=
\int_{\tilde gg\ga}i(T\ph_{\tilde g}.\tilde X_p)\dots i(T\ph_{\tilde g}.\tilde X_1)\om \\
=\int_{g(\ga)}\ph_g^*\left(i(T\ph_{\tilde g}\tilde X_p)\dots 
  i(T\ph_{\tilde g}.\tilde X_1))\om\right) 
= \int_{g\ga}i(\tilde X_p)\dots i(\tilde X_1)\ph_g^*\om =f_\ga(\ph_g^*\om).
\end{multline*}
\end{proof}  

\begin{lemma}\label{differential} 
$d_G\o f_\ga=f_\ga\o d$, where $d_G$ is the 
exterior derivative in $\Om(G)$. 
\end{lemma}

\begin{proof} Let $\om\in\Om^m(M)$. If $m<q$, by definition we have 
$d_G\o f_\ga(\om)=f_\ga\o d(\om)=0$. 

Let $\om\in\Om^{p+q}(M)$, where $p\ge 0$, and $X_1,\dots,X_p\in T_e(G)$. 
Then we have 
\begin{multline}\label{d_G}
(d_Gf_\ga(\om))(X_1^r,\dots,X_p^r)(g)=
\sum_{i=1}^p(-1)^{i-1}X_i^r(g)f_\ga(\om)(X_1^r,\dots,\widehat X_i^r,\dots,X_p^r)\\
+\sum_{i<j}(-1)^{i+j}f_\ga(\om)([X_i^r,X_j^r],X_1^r,\dots,\widehat X_i^r,\dots, 
\widehat X_j^r\dots,\dots,X_p^r)(g)\\
=\sum_{i=1}^p(-1)^{i-1}X_i^r(g)\int_\ga g^*\left(i(\tilde X_p)\dots\widehat{i(\tilde X_i)}
\dots i(\tilde X_1)\om\right)\\
+\sum_{i<j}(-1)^{i+j}\int_\ga g^*\left(i(\tilde X_p)\dots\widehat{i(\tilde X_j)}\dots
\widehat{i(\tilde X_i)}\dots i(\tilde X_1)i([\tilde X_i,\tilde X_j])\om\right)\\
=\sum_{i=1}^p(-1)^{i-1}\int_\ga g^*\left(\on{L}_{\tilde X_i}(i(\tilde X_p)\dots
\widehat{i(\tilde X_i)}\dots i(\tilde X_1)\om\right)\\
+\sum_{i<j}(-1)^{i+j}\int_\ga g^*\left( i(\tilde X_p)\dots\widehat{i(\tilde X_j)}\dots
\widehat{i(\tilde X_i)}\dots i(\tilde X_1)i([\tilde X_i,\tilde X_j])\om\right),
\end{multline}
where $\on{L}_X$ denote the Lie derivative with respect to a vector field
$X$ and, as usual, $\widehat{i(\tilde X)}$ means that the term 
$i(\tilde X)$ is omitted.

Using the formula $[\on{L}_X,i(Y)]=i([X,Y])$ it is easy to check by
induction over $p$ that for any manifold $M$ and vector fields
$X_1,\dots,X_p$ on $M$ the following formula is true.
\begin{multline*}
\sum_{i=1}^p(-1)^{i-1}\on{L}_{X_i}i(X_p)\dots\widehat{i(X_i)}\dots i(X_1)\\+
\sum_{i<j}(-1)^{i+j}i(X_p)\dots\widehat{i(X_j)}\dots\widehat{i(X_i)}\dots
i(X_1)i([X_i,X_j]\\
=i(X_q)\dots i(X_1)d+(-1)^{p-1}di(X_q)\dots i(X_1).
\end{multline*}

Applying this formula in (\ref{d_G}) we get 
$$
(d_Gf_\ga(\om_1))(X_1^r,\dots,X_p^r)(g)=\int_\ga g^*i(\tilde X_p)\dots
i(\tilde X_1)d\om=
f_\ga(d\om)(X_1^r,\dots,X_p^r)(g).
$$
\end{proof}

Consider the double complex $(C^*(G,\Om(G)),\de'_G,\de''_G)$ for the action
of the group $G$ on $G$ by left translations. Define the map
$F_\ga:(C^{**}(M,\Om(M))(C^{**}(G,\Om(G))$ as follows: for a cochain
$c\in C^p(G,\Om^q(M))$ put $F_\ga(c)=f_\ga\o c$. 

\begin{lemma}\label{d'}  
$\de'_G\o F_\ga=F_\ga\o\de'$.
\end{lemma}

\begin{proof} 
Let $c\in C^s(G,\Om^{p+q}(M)$ and $g,g_1,\dots,g_{s+1}\in G$.  
By definition we have 
\begin{multline}\label{de'_G}
(\de'_G\o F_\ga)(c))(g_1,\dots,g_{s+1})(g)=
F_\ga(c)(g_2,\dots,g_{s+1})(g)\\
+\sum_{i=1}^s(-1)^iF_\ga(c)(g_1,\dots,g_ig_{i+1},\dots, g_{s+1})(g)+
(-1)^{s+1}\on{L}_{g_{s+1}}^*F_\ga(c)(g_1,\dots,g_s)(g).
\end{multline}
For $X_1,\dots,X_p\in T_e(G)$ and $g\in G$ by (\ref{1}) and (\ref{2}) we get 
\begin{multline}\label{3}
\on{L}_{g_{s+1}}^*F_\ga(c)(g_1,\dots,g_s)(X_1^r,\dots,X_p^s)(g)\\=
\int_\ga(g_{s+1}g)^*\left(i(\widetilde{\on{ad}g_{s+1}(X_p)})\dots 
i(\widetilde{\on{ad}g_{s+1}(X_1))}c(g_1,\dots g_s)\right)\\=
\int_\ga g^*\left(i(\tilde X_p)\dots i(\tilde X_1)c(g_1,\dots,g_s)\right)=
F_\ga(g^*_{p+1}c(g_1,\dots,g_s))(X_1^r,\dots,X_p^s)(g). 
\end{multline}
Replacing the last summand in (\ref{de'_G}) by formula (\ref{3}) we get 
$$
(\de'_G\o F_\ga)(c)(g_1,\dots,g_{s+1})(g)=(F_\ga\o\de')(c)(g_1,\dots,g_{s+1})(g). 
$$
\end{proof} 

Lemmas \ref{equivariant},\ref{differential}, 
and (\ref{de'_G}) imply the following 

\begin{theorem}\label{F_} 
The map $F_\ga:C^*(G,\Om(M))\to C^*(G,\Om(G))$ is a homomorphism 
of double complexes decreasing the second degree to $q$.
\end{theorem}

Suppose that the conditions of theorem \ref{group cocycle} for an exact
$m$-form $\om$ are satisfied.
Let $\al$ be a singular smooth cycle of $M$ of dimension $m-p-1$ whose
homology class $a$
is invariant under the natural action of the group $G$ on $H_{m-p-1}(M)$.
Consider the sequence of cochains $\ph_{i,m-i-1}$ $(i=0,\dots,p)$
constructed in the proof of theorem \ref{group cocycle}. By theorem
\ref{F_} $F_\al\o\om$ is a left invariant $(p+1)$-form on $G$, i.e.,
a $(p+1)$-cocycle of the complex $C^*(\g,\mathbb R)$. Moreover, we have 
$$
\de'_GF_\al\o\ph_{i-1,m-i}+\de''_GF_\al\o\ph_{i,m-i-1}=0\quad(i=1,\dots,p).
$$
Since 
$d_G(\de'_G\o F_\al\o\ph_{p,m-p-1})=0$, for any $g\in G$ we have 
$$
c_a(\om)=\int_\ga\de'\ph_{p,m-p-1}=(F_\al\o\de'\o\ph_{p,m-p-1})(e)=
(\de'\o F_\al\o\ph_{p,m-p-1})(e).
$$

Consider the complex $(C^*(G,\mathbb R),D)$ of nonhomogeneous cochains on the 
group $G$ with values in the trivial $G$-module $\mathbb R$. Define a cochain 
$b\in C^p(G,\mathbb R)$ as follows: 
$$
b(g_1,\dots,g_p)=\int_\al\ph_{p,m-p-1}(g_1,\dots,g_p).
$$ 
By the definitions of the cocycle $c_a(\om)$ and the map $f_\ga$ and by
(\ref{de'}) we have
\begin{equation}\label{cocycle1}
\begin{split}
c_a(g_1,\dots,g_p,g)=(-1)^{p+1}&\left(f_\al\o\ph_{p,m-p-1}(g_1,\dots,g_p)(g)
-b(g_1,\dots,g_p)\right)\\&+(Db)(g_1,\dots,g_p,g).
\end{split}
\end{equation}

\begin{proposition}\label{triviality} 
Let the cycle $\ga$ be $G$-invariant. Then the cocycle $c_a(\om)$ is trivial. 
\end{proposition}
\begin{proof} By assumption we have 
$$
(F_\al\o\ph_{p,m-p-1})(g_1,\dots,g_p)(g)=\int_{g\al}\ph_{p,m-p-1}(g_1,\dots,g_p)=
b(g_1,\dots,g_p)
$$
Then by (\ref{cocycle1}) we have $c_a(\om)=Db$.
\end{proof}
 
Denote by $H$ the subgroup of $G$ consisting of all elements $g\in G$
preserving the cycle $\ga$. 
Consider the natural action of the group $G$ on the homogeneous space 
$G/H$. The projection $p_H:G\to G/H$ induces a homomorphism 
$\tilde p_H:C^*(G,\Om(G/H))\to C^*(g,\Om(G))$ of double complexes.
 
\begin{proposition}\label{G/H} 
There is a unique homomorphism of double complexes
$F_{\ga,H}:C^*(G,\Om(M))\to C^*(G,\Om(G/H))$ such that $F_\ga=\tilde p_H\o
F_{\ga,H}$.
\end{proposition}

\begin{proof} 
Note that formula (\ref{f}) implies that the form $f_\ga(\om)$
is $H$-invariant.
Moreover, by assumption for each $X\in T_e(H)$ the fundamental vector field
$\tilde X$
preserves the cycle $\ga$. 
Thus the form $i(X)\om$ vanishes on the cycle $\ga$. 
This implies that $f_\ga(\om)(X_1^r,\dots,X_p^r)=0$ whenever one of the vectors 
$X_1,\dots,X_p\in T_e(G)$ belongs to $T_e(H)$. Thus the form $f_\ga(\om)$ lies in 
the image of the map $p^*:\Om(G/H)\to\Om(G)$.  
\end{proof}
We point out the following sufficient condition of nontriviality of the cocycle
$c_a(\om)$.
\begin{theorem} \label{sufficient condition}
Let $\om\in\Om(M)^G$ be an exact $m$-form such that the conditions of theorem 
\ref{group cocycle} are satisfied and let $\al$ be a singular smooth
$(m-p-1)$-cycle on $M$ whose homology class $a$ is $G$-invariant. 
If the cohomology class of the closed left invariant form $f_\al\o\om$ in
the complex
$\Om(G)^G$ is nontrivial, the cocycle $c_a(\om)$ of the complex
$C^*(G,\mathbb R)$ is nontrivial as well.
\end{theorem}

\begin{proof} 
It is easy to check that the form $f_a\o\om$ satisfies the conditions of
proposition \ref{nontriviality} and the cocycle $c_a(\om)$ equals the
cocycle
$\de'_G\o f_\al\o\ph_{p,m-p-1}(e)$. Thus the cocycle $c_a(\om)$ is nontrivial.  
\end{proof}

Let $H$ be a subgroup $G$ preserving the cycle $\ga$. By proposition
\ref{G/H} the condition of theorem \ref{sufficient condition} can be
satisfied only if
$\dim G/H\ge m$.

\section{Continuous and differentiable cocycles}\label{continuous cocycles} 

Let $G$ be a topological group (or a Lie group which may be infinite-dimensional), 
$E$ a Frech\'et space, and $\rho:G\to\on{GL}(E)$ a representation of $G$ in $E$.
A cochain $f\in C^p(G,E)$ is continuous (differentiable) if it is a
continuous (differentiable of class $C^\infty$) map from $G^p$ to $E$. Let
$C^p_{\on{c}}(G,E)$ and $C^p_{\on{diff}}(G,E)$ denote the subspaces of continuous
and differentiable cochains of the space $C^p(G,E)$, respectively. Denote by
$H^*_{\on{c}}(G,E)$ and $H^*_{\on{diff}}(G,E)$ the cohomology of the complex
$C^*_{\on{c}}(G,E)=(C^p_{\on{c}}(G,E),\de')_{p\ge 0}$ and of
$C^*_{\on{c}}(G,E)=(C^p_{\on{diff}}(G,E),\de')_{p\ge 0}$, respectively. 
It is known
(see \cite{vE} and \cite{Gui}) that the inclusion 
$C^*_{\on{diff}}(G,E)\subset C^*_{\on{c}}(G,E)$ induces an isomorphism
$H^*_{\on{c}}(G,E)=H^*_{\on{diff}}(G,E)$ whenever $G$ is a finite
dimensional Lie group. 

Later we apply these notions to $\Om(M)$ as a
topological vector space with the $C^\infty$-topology. 
Evidently both 
$C^{**}_{\on{c}}(G,\Om(M))=(C^p_{\on{c}}(G,\Om^q(M)))$ and   
$C^{**}_{\on{diff}}(G,\Om(M))=(C^p_{\on{diff}}(G,\Om^q(M)))$  
are subcomplexes of $C^{**}(G,\Om(M))$. 

Let the conditions of theorem \ref{group cocycle} be satisfied for an exact
$m$-form $\om\in\Om(M)^G$. Assume that $G$ is a topological group or a Lie
group.
We investigate whether we can construct a sequence 
$\ph_{i,m-i-1}$ for $i=1,\dots,p$ as above which consists of continuous or
differentiable cochains. For such a sequence $c_a(\om)$ is a continuous or
differentiable cocycle.

\begin{theorem}\label{split} 
Let $M$ be a connected manifold with a countable base. Then 
for each $p>0$ we have the following decomposition in the category of topological 
vector spaces 
$$
\Om^p(M)=d\Om^{p-1}(M)\oplus H^p(M)\oplus\Om^p(M)/Z^p(M),
$$
where $Z^p(M)$ is the space of closed $p$-forms. If $H^p(M)=0$,
$d\Om^{p-1}(M)=Z^p(M)$ and $\Om^p(M)/Z^p(M)$ are Frech\'et spaces.
\end{theorem}

\begin{proof} For compact $M$ the statement follows from the Hodge decomposition 
for the identity operator 1 on $\Om^p(M)$
$1=d\o\de\o G\oplus H^p(M)\oplus\de\o d\o G$ (see, for example, \cite{deR}). 

For noncompact $M$ the statement follows from Palamodov's theorem (see \cite{Pa}, 
Proposition 5.4). 
\end{proof}

\begin{corollary} \label{continuous (differentiable) cocycle} 
Let the conditions of theorem  
\ref{group cocycle} be satisfied for an exact $m$-form $\om\in\Om(M)^G$ and
let $G$ be a topological group (a Lie group). Then one can construct a sequence 
$\ph_{i,m-i-1}$ for $i=1,\dots,p$ consisting of continuous (differentiable)
cochains and thus for a singular smooth $(m-p-1)$-dimensional cycle $\al$ whose 
homology class $a$ is $G$-invariant 
the corresponding cocycle $c_a(\om)$ is continuous (differentiable).
\end{corollary}

\begin{proof} 
The sequence $\ph_{i,m-i-1}$ $(i=1,\dots,p)$ is constructed successively  
by means of the equation $\de'\ph_{i-1,m-i}+(-1)^id\ph_{i,m-i-1}=0$. By
theorem \ref{split} for each of the above cases this equation has a
continuous (differentiable) solution
$\ph_{i,m-i-1}=L_{m-i}\o\de'\ph_{i_1,m-i}$ whenever the cochain
$\ph_{i_1,m-i}$ is continuous (differentiable).
\end{proof}

\section{Conditions of triviality of a differentiable cocycle $c_a(\om)$}

In this section we study the conditions of triviality of the cocycle
$c_a(\om)$ in the complex $C^*_{\on{diff}}(G,\mathbb R)$ whenever
$G$ is a Lie group and for the exact $m$-form $\om$ we choose a sequence of
cochains $\ph_{i,m-i-1}$ $(i=1,\dots,p)$ consisting of differentiable
cochains.
 
\begin{theorem}\label{triviality condition}
Let $M$ be a $G$-manifold, where $G$ is a Lie group preserving an exact
$m$-form $\om$, let the conditions of theorem \ref{group cocycle} be
satisfied, and for $i=1,\dots,p$
let $\ph_{i,m-i-1}\in C^i_{\on{diff}}(G,\Om^{m-i-1}(M))$. 
Then, if the cocycle $c_a(\om)$ 
is trivial, there is a cochain $f\in C^{p-1}_{\on{diff}}(G,\Om^0(G))$ such that 
$\de'_G(f_\ga\o\ph_{p-1,m-p}-d_Gf)=0$. If the group $G$ is connected, this
condition implies the triviality of the cocycle $c_a(\om)$.
\end{theorem}

\begin{proof}
Let the cocycle $c_a(\om)$ be trivial. By (\ref{cocycle1}) there is 
a cochain $\bar f\in C^p_{\on{diff}}(G,\mathbb R)$ such that for any 
$g,g_1,\dots,g_p\in G$ we have 
\begin{equation}\label{4}
(-1)^{p+1}\left(f_a\o\ph_{p,m-p-1}(g_1,\dots,g_p)(g)-b(g_1,\dots,g_p)\right)=
D\bar f(g_1,\dots,g_p,g).
\end{equation}
Define a cochain $f\in C^{p-1}_{\on{diff}}(G,\Om^0(G))$ as 
follows 
$$
f(g_1,\dots,g_{p-1})(g)=\bar f(g_1,\dots,g_{p-1},g).
$$ 
By lemma \ref{differential} we have 
\begin{multline}\label{5}
d_G((-1)^{p+1}\left((f_a\o\ph_{p,m-p-1})(g_1,\dots,g_p)(g)-b(g_1,\dots,g_p)\right)\\=
(-1)^{p+1}(f_\al\o d\ph_{p,m-p-1})(g_1,\dots,g_p)(g)\\=
\de'_G(f_\al\o\ph_{p-1,m-p})(g_1,\dots,g_p)(g).
\end{multline}

On the other hand, it is easy to check that 
\begin{equation}\label{6}
D\bar f(g_1,\dots,g_p,g)=(\de'_Gf)(g_1,\dots,g_p)(g)+(-1)^{p+1}\bar f(g_1,\dots,g_p).
\end{equation}
By (\ref{5}) and (\ref{6}), equation (\ref{4}) implies 
\begin{multline*}
\de'_G(f_\al\o\ph_{p-1,m-p})(g_1,\dots,g_p)(g)-(d_G\de'_Gf)(g_1,\dots,g_p)(g))\\=
\de'_G\left((f_\al\o\ph_{p-1,m-p}-d_Gf)\right)(g_1,\dots,g_p)(g))=0. 
\end{multline*}

Now suppose that the condition of the theorem is satisfied. We may assume that 
for any $g_1,\dots,g_p\in G$ we have $f(g_1,\dots,g_p)(e)=0$. The above
condition is equivalent to the following one
$$
\de''_G\left((-1)^{p+1}f_\al\o\ph_{p,m-p-1}-\de'_Gf\right)(g_1,\dots,g_p)(g)=0.
$$
Since the group $G$ is connected and $f(g_1,\dots,g_{p-1})(e)=0$ we have 
\begin{multline*}
\left((-1)^{p+1}f_\al\o\ph_{p,m-p-1}-\de'_Gf\right)(g_1,\dots,g_p)(g)\\=
\left((-1)^{p+1}f_\al\o\ph_{p,m-p-1}-\de'_Gf\right)(g_1,\dots,g_p)(e)\\=
(-1)^{p+1}(b-\bar f)(g_1,\dots,g_p).
\end{multline*}
Using (\ref{cocycle1}) and (\ref{6}) we get 
$(c_a(\om)-D\bar f)(g_1,\dots,g_p,g)=0$. 
This concludes the proof.  
\end{proof}  
    
\begin{corollary}\label{triviality p=1} 
Let the conditions of theorem \ref{triviality condition} be satisfied 
for $m=2$ and $p=1$. Then, if the 2-cocycle $c_a(\om)$ is trivial, there is
a smooth function $f$ on $G$ such that the 1-form $f_\al\o(\ph_{0,1})-d_Gf$
on $G$ is left invariant. In particular, the cohomology class of the form
$f_\al\o\om$ in the complex $\Om(G)^G$ is trivial. If the group $G$ is
connected, the above condition implies the triviality of the cocycle
$c_a(\om)$.
\end{corollary} 

\begin{example}\label{R^n}
Consider the abelian group $G=\mathbb R^n$ acting on itself by translations.
Evidently an $m$-form 
$$
\om(x)=\sum_{i_1<\dots<i_m}\om_{i_1\dots i_m}(x)dx_{i_1}\wedge\dots\wedge dx_{i_m},
$$
where $x=(x_1,\dots,x_n)\in\mathbb R^n$, is $G$-invariant if and only if the
coefficients
$\om_{i_1\dots i_p}$ are constant. Then the differential of the complex
$\Om(\mathbb R^n)^G$ is trivial.

Evidently the conditions of theorem \ref{group cocycle} are satisfied for
each nonzero $m$-form $\om$ with constant coefficients on $\mathbb R^n$ for
$p=m-1$.
It is easy to check that the sequence  of 
cochains
$\ph_{j,m-j-1}\in C^j_{\on{diff}}(G,\Om^{m-j-1}(\mathbb R^n))$
$(j=0,\dots,m-1)$ corresponding to $\om$ 
can be defined as follows.
$$
\ph_{j,m-j-1}(a_1,\dots,a_j)=
\frac{(-1)^{j-1}}{m(m-1)\dots(m-j)}i(x)i(a_j)i(a_{j-1})\dots i(a_1)\om,
$$
where $a_k=(a_{k,1},\dots,a_{k,n})\in\mathbb R^n$ $(k=1,\dots,j)$ and on
the right hand side we consider each $a_k$ as a constant vector field on
$\mathbb R^n$ and $x$ as identical vector field on $\mathbb R^n$. Then we
have
$$
\de'\ph_{m-1,0}(a_1,\dots,a_m)=\frac1{m!}i(a_m)i(a_{m-1})\dots i(a_1)\om,
$$ 
where $a_1,\dots,a_m\in\mathbb R^n$.

Take the point $0\in\mathbb R^n$ as the cycle $\al$. Then 
$F_\al:\Om(\mathbb R^n)^G\to\Om(G)$ is the identity map. By proposition 
\ref{sufficient condition} 
the cocycle $c_a(\om)$ is nontrivial in the complex $C^*(G,\mathbb R)$. 
\end{example}
  
\section {Cocycles on groups of diffeomorphisms}
In this section we indicate nontrival cocycles for groups of
diffeomorphisms of a manifold preserving a family of exact forms.

Let $(\om_i)_{i\in I}$ be a family of smooth differential forms on a
manifold $M$. Denote by $\on{Diff}(M,(\om_i))$ the group of diffeomorphisms
of $M$ preserving all forms $\om_i$. We consider $\on{Diff}(M,(\om_i))$ as
a topological group with respect to $C^\infty$-topology or as a
infinite-dimensional Lie group if such a structure on
$\on{Diff}(M,(\om_i))$ exists. By proposition \ref{continuous
(differentiable) cocycle} one can suppose that the cocycle $c_a(\om_i)$ is
a continuous or differentiable cocycle.

Let $(\om_i)_{i\in I}$ be a family of nonzero  differential forms on
$\mathbb R^n$ with constant coefficients. Put $G=\on{Diff}(\mathbb
R^n,(\om_i))$.

\begin{theorem}\label{Diff(R^n,om)} 
Let $(\om_i)_{i\in I}$ be a family of nonzero differential forms on 
$\mathbb R^n$ with constant coefficients such that
$\deg \om_i=d_i$.  
Then, for each $i\in I$ and for $0\in\mathbb R^n$ as a zero dimensional
cycle $\al$, the cocycle $c_a(\om_i)$ is defined and nontrivial in the
complex $C^*(G,\mathbb R)$.
\end{theorem}    

\begin{proof} 
Evidently the conditions of theorem \ref{group cocycle} are satisfied for
each form $\om_i$ and $p=m_i-1$ and then the cocycle $c_a(\om_i)$ of the
complex
$C^*(G,\mathbb R)$ is defined.  Obviously the group $G$ contains
the group $\mathbb R^n$ acting on $\mathbb R^n$ by translations. 
Consider the restriction of the cocycle $c_a(\om)$ to the
subgroup
$\mathbb R^n$. By example \ref{R^n} this restriction is nontrivial. Then
the cocycle
$c_a(\om_i)$ is nontrivial as well. 
\end{proof}

\begin{corollary}\label{R^n x M} Let $M$ be a connected manifold such that 
$$
H^1(M,\mathbb R)=\dots=H^{m-1}(M,\mathbb R)=0,
$$ 
let $\om_0$ be a nonzero $m$-form on $\mathbb R^n$ with constant
coefficients, and let $\om_M$ be an exact $m$-form on $M$.
Consider $\om=\om_0+\om_M$ as an $m$-form on $\mathbb R^n\x M$. Then for
the group
$G=\on{Diff}(\mathbb R^n\x M,\om)$ and a point $x\in\mathbb R^n\x M$ as a
zero dimensional cycle $\al$, the cocycle $c_a(\om)$ is defined and
nontrivial in the complex $C^*(R^n\x M,\mathbb R))$. 
\end{corollary}

\begin{proof} Evidently the conditions of theorem \ref{group cocycle} for
the form $\om$ are satisfied and then the cocycle $c_a(\om)$ of the complex
$C^*(G,\mathbb R)$ is defined.
Consider the group $\on{Diff}(\mathbb R^n,\om_0)$ acting on the first
factor of
$\mathbb R^n\x M$ as a subgroup of $G$ and the restriction of the cocycle
$c_a(\om)$ to this subgroup. Since
the subgroup $\on{Diff}(\mathbb R^n,\om_0)$ preserves the form $\om_M$ as a
form on
$\mathbb R^n\x M$ by theorem \ref{Diff(R^n,om)}, this restriction is a
nontrivial cocycle. Thus $c_a(\om)$ is a nontrivial cocycle of the complex
$C^*(G,\mathbb R)$ as well. 
\end{proof}

We indicate the following partial case of corollary \ref{R^n x M} (see also
\cite{LM}).

\begin{corollary} Let $\om_0$ be the standard symplectic 2-form on $\mathbb
R^{2n}$ and
let $m=1,\dots,n$. Let $M$ be a connected manifold such that $2m\le\dim M$,
$$
H^1(M,\mathbb R)=\dots=H^{2m-1}(M,\mathbb R)=0,
$$ 
and $\om_M$ an exact $2m$-form on a manifold $M$. Consider the $2m$-form
$\om=\om_0^m+\om_M$ on $\mathbb R^{2n}\x M$. Then for the group
$G=\on{Diff}(\mathbb R^{2n}\x M,\om)$ and a point $x\in\mathbb R^{2n}\x M$
as a zero dimensional cycle $\al$, the cocycle $c_a(\om)$ is defined and
nontrivial in the complex $C^*(R^{2n}\x M,\mathbb R))$.
\end{corollary} 

Let $M$ be a connected compact oriented manifold with a volume form
$\upsilon_M$ such that $\int_M\upsilon_M=1$.

\begin{theorem}\label{M,upsilon} Let $(\om_i)_{i\in I}$ be a family of
nonzero differential forms on $\mathbb R^n$ with constant coefficients such
that $\deg \om_i=d_i$. Consider the family 
$\{\om_i\wedge\upsilon_M\}_{i\in I}$ of forms on $\mathbb R^n\x M$. 
For a cycle $\al=0\x M$ 
of the homology class $a$ on $\mathbb R^n\x M$ and each $i\in I$ the
cocycle $c_a(\om_i\wedge\upsilon_M)$ on
the group $G=\on{Diff}(\mathbb R^n\x M,(\om_i\wedge\upsilon_M))$ with
values in the trivial $G$-module $\mathbb R$ is defined and nontrivial.
\end{theorem}

\begin{proof}  
Evidently the conditions of theorem \ref{group cocycle} for
each form
$\om_i\wedge\upsilon_M$ and $p=n-1$ are satisfied. Then the cocycle 
$c_a(\om_i\wedge\upsilon_M)$ of the complex $C^*(G,\mathbb R)$ is defined. 

Consider the group $\on{Diff}(\mathbb R^n,(\om_i))$ as a subgroup of 
the group $G$ and the restriction of the cocycle $c_a(\om)$ to this
subgroup. Since the
subgroup $\on{Diff}(\mathbb R^n,(\om_i))$ preserves the form $\upsilon$ as
a form on
$\mathbb R^n\x M$, as in the proof of corollary \ref{R^n x M}
we can show that the
above restriction is a nontrivial cocycle. Thus the cocycle $c_a(\om)$ is 
nontrivial as well.
\end{proof}
 
We indicate the following partial case of theorem \ref{M,upsilon}.

\begin{corollary}\label{volume} 
Let $\upsilon_0$ be the standard volume form on $\mathbb R^n$.
Then for the cycle $a=0\x M$ the cocycle $c_a(\upsilon_0\wedge\upsilon_M)$
on the group $G=\on{Diff}(\mathbb R^n\x M,\upsilon_0\wedge\upsilon_M)$ is
defined and nontrivial in the complex
$C^*(G,\upsilon_0\wedge\upsilon_M),\mathbb R)$.
\end{corollary}

We consider the space $\mathbb C^{2n}$ and a skew-symmetric bilinear form
$\om$ of rank 2 on it. Let $\tilde\om$ be the differential
2-form corresponding to $\om$ on
$\mathbb C^{2n}$ as a complex manifold. By definition the form $\tilde\om$
has constant
coefficients.  Then $\tilde\om=\tilde\om_1+i\tilde\om_2$, where
$i=\sqrt -1$, $\tilde\om_1$ and $\tilde\om_2$ are real differential 2-forms
with constant coefficients on $\mathbb C^{2n}=\mathbb R^{4n}$ as a real
$4n$-dimensional manifold.
 
Similarly, consider the space $\mathbb C^n$, a sqew-symmetric $n$-form
$\upsilon$ of maximal rank on it, and the corresponding differential
$n$-form $\tilde\upsilon$ on
$\mathbb C^n$. Consider the real differentiable $n$-forms $\tilde\upsilon_1$
and $\tilde\upsilon_2$ on $\mathbb C^n=\mathbb R^{2n}$
as a real $2n$-dimensional manifold defined by an equality
$\tilde\upsilon=\tilde\upsilon_1+i\tilde\upsilon_2$.

In both cases above we can apply theorem \ref{Diff(R^n,om)} to the forms $\tilde\om_1$,  $\tilde\om_2$, $\tilde\upsilon_1$ and $\tilde\upsilon_2$. We leave to the reader to define the 
corresponding cocycles on the group of diffeomorphisms preserving the forms 
$\tilde\om$ and $\tilde\upsilon$ and to formulate the statement similar to those of 
corollaries \ref{R^n x M} and \ref{M,upsilon}.

\end{document}